\documentclass[10pt, titlepage]{amsart}

\usepackage{ae} 
\usepackage[T1]{fontenc}
\usepackage[cp1250]{inputenc}
\usepackage{amsmath}
\usepackage{amssymb, amsfonts,amscd,verbatim}
\usepackage[dvips]{graphicx}
\usepackage{latexsym}
\usepackage{indentfirst}
\usepackage{latexsym}

\usepackage{graphicx}
\usepackage{verbatim}   
\usepackage{color}      
\usepackage{subfigure}  

\usepackage{amsmath}    

\theoremstyle{plain}
\newtheorem{Prop}{Proposition}[section]
\newtheorem{Thm}[Prop]{Theorem}

\theoremstyle{definition}
\newtheorem{Def}[Prop]{Definition}

\theoremstyle{remark}

\def\int{\mathop{\roman{int}}}

\def\1{^{-1}}

\def\diam{\text{diam}}

\def\asdim{\mathrm{asdim}}

\def\diam{\mathrm{diam}}
\def\dokaz{{\bf Proof. }}
\def\edokaz{\hfill $\blacksquare$}

\numberwithin{equation}{section}

\input pstricks.tex
\input xy
\xyoption{all}

\begin{document}
\title[
Property A and asymptotic dimension
]%
   {Property A and asymptotic dimension}

\author{M.~Cencelj}
\address{IMFM,
Univerza v Ljubljani,
Jadranska ulica 19,
SI-1111 Ljubljana,
Slovenija }
\email{matija.cencelj@guest.arnes.si}

\author{J.~Dydak}
\address{University of Tennessee, Knoxville, TN 37996, USA}
\email{dydak@math.utk.edu}

\author{A.~Vavpeti\v c}
\address{Fakulteta za Matematiko in Fiziko,
Univerza v Ljubljani,
Jadranska ulica 19,
SI-1111 Ljubljana,
Slovenija }

\email{ales.vavpetic@fmf.uni-lj.si}

\date{ \today
}


\thanks{Supported in part by the Slovenian-USA research grant BI--US/05-06/002 and the ARRS
research project No. J1--6128--0101--04}
\thanks{The second-named author was partially supported
by MEC, MTM2006-0825.}

\begin{abstract}

The purpose of this note is to characterize the asymptotic dimension
$\asdim(X)$
of metric spaces $X$ in terms similar to Property A of Yu \cite{Yu00}:

\begin{Thm}
If $(X,d)$ is a metric space and $n\ge 0$, then the following conditions are equivalent:
\begin{itemize}
\item[a.] $\asdim(X,d)\leq n$,
\item[b.] For each $R,\epsilon > 0$ there is $S > 0$
and finite non-empty subsets $A_x\subset B(x,S)\times N$, $x\in X$,
such that $\frac{\vert  A_x\Delta A_y\vert}{\vert A_x\cap A_y\vert} < \epsilon$ if $d(x,y) < R$ and the projection
of $A_x$ onto $X$ contains at most $n+1$ elements
for all $x\in X$,
\item[c.] For each $R > 0$ there is $S > 0$
and finite non-empty subsets $A_x\subset B(x,S)\times N$, $x\in X$,
such that $\frac{\vert  A_x\Delta A_y\vert}{\vert A_x\cap A_y\vert} < \frac{1}{n+1}$ if $d(x,y) < R$ and the projection
of $A_x$ onto $X$ contains at most $n+1$ elements
for all $x\in X$.

\end{itemize}

\end{Thm}

\end{abstract}

\maketitle

\medskip
\medskip
\tableofcontents

\section{Introduction}
Property A was introduced by G.Yu in \cite{Yu00}
in order to prove a special case of the Novikov Conjecture. We adopt the following definition from \cite{NY} (see also \cite{W}):
\begin{Def}
A discrete metric space $(X, d)$ has
property A if for all $R, \epsilon > 0$, there exists a family $\{A_x\}_{x\in X}$ of finite, non-empty
subsets of $X \times N$ such that:
\begin{itemize}
\item for all $x,y\in X$ with $d(x,y)\leq R$
we have
$\frac{\vert  A_x\Delta A_y\vert}{\vert A_x\cap A_y\vert} < \epsilon$
\item there exists $S > 0$ such that for each $x\in X$,
if $(y,n)\in A_x$, then $d(x,y)\leq S$
\end{itemize}

\end{Def}

Asymptotic dimension was introduced by M. Gromov
in \cite{Grom} (see section 1.E) as a large-scale analogue of the classical notion of topological covering dimension. It is a coarse invariant that has been extensively investigated (see chapter 9 of 
\cite{Roe lectures} for some results
and further references).

\begin{Def}
A metric space $X$ is said to have finite asymptotic dimension
if there exists $k \ge 0$ such that for all $L > 0$ there exists a uniformly
bounded cover of $X$ (that means the existence of $S > 0$ such that all elements of the cover are of diameter at most $S$) of Lebesgue number at least $L$ (that means every $R$- ball $B(x,R)$ is contained in some element of the cover) and multiplicity at most $k +1$ (i.e. each point of $X$ belongs to at most $k+1$ elements of the cover). The
least possible such $k$ is the {\bf asymptotic dimension} of $X$.
\end{Def}

One of the basic results is that spaces of finite asymptotic dimension
have property A and known proofs of it use Higson-Roe
characterization of Property A (see \cite{HR} and \cite{W}).
The purpose of this note is to provide a simple proof of that result
and prove a characterization of asymptotic dimension in terms similar to Property A.

\section{Main theorem}

\begin{Thm}\label{MainThm}
If $(X,d)$ is a metric space and $n\ge 0$, then the following conditions are equivalent:
\begin{itemize}
\item[a.] $\asdim(X,d)\leq n$,
\item[b.] For each $R,\epsilon > 0$ there is $S > 0$
and finite non-empty subsets $A_x\subset B(x,S)\times N$, $x\in X$,
such that $\frac{\vert  A_x\Delta A_y\vert}{\vert A_x\cap A_y\vert} < \epsilon$ if $d(x,y) < R$ and the projection
of $A_x$ onto $X$ contains at most $n+1$ elements
for all $x\in X$,
\item[c.] For each $R > 0$ there is $S > 0$
and finite non-empty subsets $A_x\subset B(x,S)\times N$, $x\in X$,
such that $\frac{\vert  A_x\Delta A_y\vert}{\vert A_x\cap A_y\vert} < \frac{1}{n+1}$ if $d(x,y) < R$ and the projection
of $A_x$ onto $X$ contains at most $n+1$ elements
for all $x\in X$.

\end{itemize}

\end{Thm}

\dokaz a)$\implies$ b).
Suppose $\asdim(X,d)\leq n$ and $R, \epsilon > 0$.
Pick a uniformly bounded cover $\mathcal{U}$ of $X$
of multiplicity at most $n+1$ and Lebegue number
at least $L=2R+\frac{2R\cdot n}{\epsilon}$.
Let $S$ be a number such that $\diam(U) < S$
for each $U\in \mathcal{U}$.
Pick $a_U\in U$ for each $U\in \mathcal{U}$
and define $A_x$ as the union of sets $a_U\times \{1,\ldots l_U(x)\}$, where $x\in U$ and $l_U(x)$ is the length
of the shortest $R$-chain joining $x$ and a point
outside of $U$ (if there is no such chain, we put $l_U(x)$
equal to the integer part of $\frac{L}{R}+1$).
If $d(x,y) < R$, then $\vert l_U(x)-l_U(y)\vert \leq 1$,
so $\vert  A_x\Delta A_y\vert \leq 2n$ (as the total number of
elements of $\mathcal{U}$ containing exactly one of $x$ or $y$
is at most $2n$),
and $\vert A_x\cap A_y\vert > \frac{L-R}{R}-1$
(choose $U$ containing $B(x,L)$ and notice every $R$-chain joining
$x$ or $y$ to $X\setminus U$ must have at least $\frac{L-R}{R}$ elements),
yielding $\frac{\vert  A_x\Delta A_y\vert}{\vert A_x\cap A_y\vert} < \frac{2n\cdot R}{L-2R} \le \epsilon$.
\par c)$\implies$ a). 
Given $R > 0$  pick $S > 0$
and finite subsets $A_x\subset B(x,S)\times N$, $x\in X$,
such that $\frac{\vert  A_x\Delta A_y\vert}{\vert A_x\cap A_y\vert} < \frac{1}{n+1}$ if $d(x,y) < R$ and the projection
of $A_x$ onto $X$ contains at most $n+1$ elements
for all $x\in X$.
Define sets $U_x$ as consisting precisely of
$y\in X$ such that $(\{x\}\times N)\cap A_y\ne \emptyset$.
The multiplicity of the cover $\{U_x\}_{x\in X}$ of $X$
is at most $n+1$ as $z\in \bigcap\limits_{i=1}^k U_{x_i}$ implies
$x_i$ belongs to the projection of $A_z$, so $k\leq n+1$.
Given $x\in X$ choose $z\in X$ so that
$\vert(\{z\}\times N)\cap A_x\vert$ maximizes
all $\vert(\{y\}\times N)\cap A_x\vert$.
In particular $\vert(\{z\}\times N)\cap A_x\vert \ge \frac{\vert A_x\vert}{n+1}$.
If $d(x,y) < R$ we claim $y\in U_z$ which proves that the Lebegue number of $\{U_x\}_{x\in X}$ is at least $R$.
Indeed, $y\notin U_z$ implies $\vert A_x\Delta A_y\vert \ge 
\frac{\vert A_x\vert}{n+1}$, so 
$\frac{\vert  A_x\Delta A_y\vert}{\vert A_x\vert} \ge \frac{1}{n+1}$,
a contradiction.
\edokaz


\begin{thebibliography}{99}


\bibitem{Grom}
M. Gromov, {\em Asymptotic invariants for infinite groups}, in
Geometric Group Theory, vol. 2, 1--295, G. Niblo and M. Roller,
eds., Cambridge University Press, 1993.

\bibitem{HR} 
N. Higson and J. Roe, {\em Amenable group actions and the Novikov conjecture},
J. Reine Agnew. Math. 519 (2000).

\bibitem{NY}
Piotr Nowak and Guoliang Yu,
{\em What is ... Property A?},
Notices of the AMS Volume 55, Number 4, pp.474--475.

\bibitem{Roe lectures}
J. Roe, {\em Lectures on coarse geometry}, University Lecture
Series, 31. American Mathematical Society, Providence, RI, 2003.

\bibitem{W}
Rufus Willett, {\em Some notes on Property A}, arXiv:math/0612492v2 [math.OA]

\bibitem{Yu00}
G. Yu, {\em The coarse Baum-Connes conjecture for spaces which admit a uniform embedding into Hilbert space}, Inventiones 139 (2000), pp. 201--240.

\end{thebibliography}
\end{document}